\documentclass[12pt, draft]{article}
\usepackage[cp1251]{inputenc}
\usepackage{amssymb,amsmath,eucal,latexsym,amsthm}
\usepackage{indentfirst}
\usepackage{enumerate}
\usepackage{verbatim}
\usepackage[dvips]{graphicx}
\usepackage[centerlast,small]{caption2}

 \theoremstyle{plain}
\newtheorem{Theorem}{Theorem}[section]

\newtheorem*{Hypothesis}{Hypothesis}

\theoremstyle{definition}

\theoremstyle{remark}

\newcommand{\B}{\mathbb{B}}
\newcommand{\D}{\mathbb D}
\newcommand{\R}{\mathbb{R}}
\newcommand{\bC}{\mathbb{C}}
\newcommand{\N}{\mathbb{N}}
\newcommand{\e}{\varepsilon}

\DeclareMathOperator{\const}{const}
\DeclareMathOperator{\dist}{dist} \DeclareMathOperator{\Mer}{Mer}

\title{The representation of  a meromorphic function
as the quotient of entire functions and\\
Paley problem in ${\bC}^n$: survey of some
results}  
\author{B.~N. Khabibullin
\thanks{Supported  by the Russian Foundation of Basic Research
under Grant No.~00-01-00770.}\\
\small{Dedicated to Professor N.~I.~Akhiezer in the year
of his 100th anniversary}}
\date{}

\begin{document}
\maketitle
\begin{abstract}
The classical representation problem for a meromorphic function $f$
in  $\bC^n$, $n\ge 1$,  consists in representing $f$ as the quotient $f=g/h$
of two entire functions $g$ and $h$, each with logarithm of
modulus majorized by a function as close as possible to the
Nevanlinna characteristic. Here we introduce  generalizations of the
Nevanlinna characteristic and give a short survey of classical and recent
results on the representation of a meromorphic function  in terms such
characteristics. When $f$ has a finite lower order, the Paley problem
on  best possible estimates of the growth of entire functions
$g$ and $h$ in the representations $f=g/h$  will be considered.
Also we point out to  some unsolved problems in this area.
\end{abstract}

\pagestyle{myheadings}
\markboth{B.~N. KHABIBULLIN}{THE REPRESENTATION OF A MEROMORPHIC
FUNCTION...}

\subsection*{Introduction}

Let $f$ be a meromorphic function  (we write
$f\in \Mer$)  in the complex plane.
At the beginning of the last century the definition of the
classic Nevanlinna characteristic  $T(r; f)$ and the Liouville Theorem
on the growth of entire functions were giving
an one of first results on the representation of a meromorphic
function as the quotient of entire functions.
\begin{Theorem}[J.~Liouville--R.~Nevanlinna]\label{LN}
A meromorphic function $f$ in the complex plane can be represented
as the quotient $f=g/h$ of polinomials  $g$ and $h$, i.e.
$$
\log \bigl( |g(z)|+|h(z)|\bigr) \leq O( \log r)\, , \quad r\to +\infty
\, ,
$$
iff
$$
T(r; f)\le O( \log r)\,  , \quad r\to +\infty \, .
$$
\end{Theorem}

The classic Nevanlinna Theorem on the representation of a meromorphic
function $f$ of bounded characteristic  $T(r; f)\le \const \,$, $r<1$,
in unit disk $\D$ as the quotient of bounded holomorphic functions in $\D$
can be considered also as one of sources of our survey.

A first essential result for transcendental meromorphic functions
follows immediately from the Lindel\"of Theorem of the beginning
of XX century (see, for example, \cite{Rubel}, \cite{GO}):

\begin{Theorem}[E.~Lindel\"of]\label{Li}
A  meromorphic function $f$ in the complex plane can be represented
as the quotient $f=g/h$ of entire functions  $g$ and $h$
of finite type for order $\rho$, i.e.
$$
\log \bigl( |g(z)|+|h(z)|\bigr) \leq O( r^{\rho})\, , \quad r\to
+\infty \, ,
$$
iff
$$
T(r; f)\le O(  r^{\rho})\,  , \quad r\to +\infty \, .
$$
\end{Theorem}

Here we consider  far-reaching generalizations of
these classic results and  give a  survey of basic and
recent results on the {\it problem of the representation\/}
of a meromorphic function in $\bC^n$
as the quotient of entire functions in terms of  various characteristics.
Also we point out to  some unsolved problems in this area.

Let $f \in \Mer$ in $\bC^n$ and
\begin{equation}
\label{repr}
f=\frac{g_f}{h_f}
\end{equation}
be a {\it canonical  representation\/} of $f$ as the quotient of
entire functions $g_f$ and $h_f$ which are locally relatively prime, i.e.
at each  point $z\in \bC^n$, where $g_f(z)=h_f(z)=0$, the functions
$g$ and $h$ are relatively prime in the ring of germs of functions analytic at
$z$. In $\bC^n$ such representations always exist \cite{GunnR},
\cite{Hormander}. Both functions
\begin{subequations}
\label{uf}
\begin{align}
u_f &=\max \bigl\{ \log |g_f|, \log |h_f|\bigr\} \label{ufi}  \\
\notag \text{ or } \qquad &\\
u_f &=\log \bigl( |g_f|^2 + |h_f|^2 \bigr)^{1/2}  \label{ufii}
\end{align}
\end{subequations}
are plurisubharmonic on $\bC^n$.
Consideration of these {\it characteristic functions\/} $u_f$ is useful and
natural, since  $u_f$ is used to define the various types
of Nevanlinna characteristic in particular in the Shimizu-Ahlfors form etc.
\cite{GO}, \cite{Hayman},  \cite{Kondratyuk}, \cite{Kujala},
\cite{Khabib93-I}, \cite{Skoda}--\cite{Skoda72}, \cite{Stoll}--\cite{Taylor}.

A function $g$ on $\bC^n$ is said to be {\it circular\/}
if $g(e^{i\theta}z)=g(z)$ for all $z\in \bC^n$ and $\theta \in \R$.

First we introduce the {\it circular Nevanlinna characteristic}
\begin{equation}
\label{Tcirc}
{T_f^c}(z)=\frac1{2\pi} \int\limits_{0}^{2\pi} u_f(e^{i\theta}z)\, d\theta \, ,
\quad z\in \bC^n\, ,
\end{equation}
which is a circular plurisubharmonic function \cite{Kujala}.
For example, if  $f(0)=1$ and $T_f(t; \zeta )$, $t\geq 0$,
are the Nevanlinna characteristics
(in the Shimizu-Ahlfors form for the case \eqref{ufii}) \cite{GO}
of the family of the meromorphic functions $f_{\zeta}(w)=f(w\zeta )$
of one variable $w\in \bC$, where $\zeta$ runs through the unit sphere
$S^n \subset \bC^n$, then $T_f(t; \zeta )\leq {T_f^c}(z)$,
$z=te^{i\theta}\zeta$. Equality holds
if $n=1$ or if $g_f$ and $h_f$ do not have common zeros ($n>1$).

Let $\mathbf{\mbox{\L}} =(L_1, \ldots , L_k)$  be a
fixed simply ordered collection of complex vector subspaces
in ${\bC}^n$ and ${\bC}^n$ is the direct sum
of $L_k$. As usual, we write
\begin{equation}
\label{directsum}
\bC^n = L_1\oplus \cdots \oplus L_k\, .
\end{equation}

Denote by $\R_+$ the set of nonnegative numbers.

Let $z\in \bC^n$  and
\begin{subequations}\label{z}
\begin{align}
z&=w_1+\cdots +w_k, \qquad  w_p\in L_p,\quad p=1, \ldots , k, \label{zz}\\
r_p&=|w_p|,\quad \vec r\, =(r_1, \ldots , r_k)\in \R_+^k, \label{zw}
\end{align}
\end{subequations}
where $|w_p|$ is the Euclidean norm $w_p$ in $L_p$.
For every $z\in \bC^n$  the representation \eqref{zz} is unique.

If $u(z)$ is a function on $\bC^n$,
then we consider   the function $u$ also as the function
$u(z)=u(w_1, \ldots , w_k)$  on the direct sum \eqref{directsum}
under the  notations \eqref{z}. We keep also the
same notation $u$ for $u(w_1, \ldots , w_k)$.

We introduce {\it Nevanlinna\/} $\mathbf{\mbox{\L}}$-{\it characteristic\/} of
$f\in \Mer$ in the form
$$
T_{\mathbf{\mbox{\L}}}(\vec{r}\, ; f)=\int\limits_{S_1} \cdots \int\limits_{S_k}
u_f(r_1\zeta_1, \ldots , r_k\zeta_k)
\, ds^{(1)} (\zeta_1)\cdots ds^{(k)} (\zeta_k),
\quad \vec{r}\, \in \R_+^k \, ,
$$
where the functions $u_f$ are defined in \eqref{uf},
and $ds^{(p)}$ is the area element on the unit sphere $S_p$ in $L_p$.

The {\it Nevanlinna\/} $\mathbf{\mbox{\L}}$-{\it characteristic\/} can
be defined also in the form
\begin{equation}
\label{T}
T_{\mathbf{\mbox{\L}}}(\vec{r}\, ; f)=\int\limits_{S_1} \cdots \int\limits_{S_k}
T_f^c(r_1\zeta_1, \ldots , r_k\zeta_k)
\, ds^{(1)} (\zeta_1)\cdots ds^{(k)} (\zeta_k).
\end{equation}
\newcommand{\Fl}[3]{#3_{#1},\ldots ,#3_{#2}}
\newcommand{\Fc}[3]{#3_{#1}\cdots #3_{#2}}

It is easy to check that every Nevanlinna
$\mathbf{\mbox{\L}}$-characteristic is defined by \eqref{repr}--\eqref{uf}
to within an additive  constant.

The Nevanlinna $\mathbf{\mbox{\L}}$-characteristic gives the various variants
of the classical Nevanlinna characteri\-s\-t\-ics of $f$ in $\bC^n$.
So, if $k=1$ and $L_1=\bC^n$ in \eqref{directsum},
then we get Nevanlinna characteristic $T(r; f)$, $r\geq 0$,
by exhaustion of $\bC^n$ by balls $r\B^n$, where $\B^n$ is the unit ball in
$\bC^n$; if $k=n$ and $L_1=\cdots =L_n=\bC$ in \eqref{directsum},
then we get Nevanlinna characteristic $T(\vec r\, ; f)$,
$\vec r\, \in \R_+^n$, by exhaustion of $\bC^n$ by  polydisks
$P^n(\vec r\, )=\{z\in \bC^n: |w_p|=r_p, p=1, \ldots , n \}$
with the preceding notations \eqref{z}
(see  \cite{Stoll}--\cite{Taylor}, \cite{Kujala}, \cite{Skoda}--\cite{Skoda72},
\cite{Khabib93-I} for various variants of definitions of such
Nevanlinna characteristics).

Below we formulate the well-known and recent main results on the
representation of a meromorphic function on $\bC^n$, $n\ge 1$.

The author expresses deep gratitude to
I.~F.~Krasichkov-Ternovski\u{\i} for useful discussion
of this review and a correction of the text, and to M.~L.~So\-din,
A.~A.~Kondratyuk and  Ya.~V.~Vasyl'kiv for a valuable information,
and also to  the reviewer for helpful remarks.

\section{The representation  of arbitrary\\ meromorphic functions}

\setcounter{equation}{0}

\renewcommand{\thefootnote}{\fnsymbol{footnote}}

\addtocounter{footnote}{+1}

A continuous increasing nonnegative function $\lambda $ on $\R_+$
is called a {\it growth function\/} (see \cite{Taylor-Rubel},
\cite{Miles} and \cite{Kujala}).

For $n=1$ and  for a {\it meromorphic function $f$ of  finite $\lambda$-type},
i.e. satisfying the condition
\begin{equation}
\label{ballestm}
T(r; f)\le \const \cdot \lambda (r)+\const \, , \quad r>0\, ,
\end{equation}
where $\lambda $ has either  the {\it slow growth\/}, i.e. satisfies
the condition
\begin{equation}
\label{slow1}
\lambda (2r)\le O\bigl(\lambda (r)\bigr) \,  ,\quad r \to +\infty\; ,
\end{equation}
or  $\log \lambda (\exp t)$ is {\it convex\/}, L.~A.~Rubel and B.~A.~Taylor
\cite[Theorems 5.4, 3.5-6, 3.2]{Taylor-Rubel} indicated (1968)
{\it a construction of entire functions $g$ and $h$ such that $f=g/h$
and
\begin{equation}\label{balleste}
\log \bigl(|g(z)|+|h(z)|\bigr)\le A\lambda (Br)+C \; ,\quad r=|z|>0\; ,
\end{equation}
where $A$, $B$ and C are constants\/}.

All conditions for $\lambda$  was removed for $n=1$ (1970-72) by
J.~B.~Miles in \cite{MilesB}, \cite{Miles}
(see also the original proof of J.~B.~Miles
in \cite[Ch.~4]{Kondratyuk}, \cite[Ch.~14]{Rubel}):

\begin{Theorem}[L.~A.~Rubel--B.~A.~Taylor--J.~B.~Miles]
\label{RTM}
 Every fu\-n\-c\-t\-i\-on
$f\in \Mer$ in the complex plane $\bC$ can be represented as the quotient
 $f=g/h$ of entire functions $g$ and $h$ such that
\begin{equation}
\label{ballesteM}
\log \bigl(|g(z)|+|h(z)|\bigr)\le AT(Br; f)+C \; ,\quad r=|z|>0\; ,
\end{equation}
where $A$, $B$ and $C$ are
constants\footnote{The importance of such
results for $n>1$ was marked by W.~Stoll in his book  \cite[p.~85]{Stoll85}:
{\it `` Unfortunately the theorem of Miles  {\rm \cite{Miles}} has not been
proved for several variables''}.}.
\end{Theorem}

Analogous results for functions in the half-plane
and in the unit disk can be find  in \cite{Mlt} and in \cite{Beck} recp..

The following special case of the representation theorems
was established for $n=1$ in \cite{Gol'dberg} (1972):
\begin{Theorem}[A.~A. Gol'dberg]
\label{Gol}
Every meromorphic function $f$ in the complex plane
${\bC}$ can be re\-p\-r\-e\-s\-e\-n\-ted
as the quot\-i\-ent $f=g_1/g_2$ of entire functions $g_1$ and $g_2$
\underline{without  common zeros} such that
$$
\log T(r; g_k)=o\bigl(T(r; f)\bigr) \, , \quad r\to \infty \; ,\quad k=1,2\; .
$$
\end{Theorem}

This theorem (for functions $g_1$ and $g_2$ {\it without common zeros})
is, in a certain sense,  best possible (see also \cite{Gol'dberg}).
Further results in this direction (\underline{without  common zeros})
can to find  in \cite{Skoda3}, \cite{Berg}, \cite{VSh}, \cite{KonVas}.

Our last one-dimensional  result \cite{Khabib01} (2001) is the following
\begin{Theorem}[B.~N. Khabibullin]
\label{Khab01}
Let $f$ be meromorphic function in the complex plane.
For each  nonincreasing convex function $\e (r)>0$ in $\R_+$
entire functions $g$ and $h$ in the representation $f=g/h$
can be chosen so that the estimate
\begin{equation}
\label{Te}
\log \bigl(|g(z)|+|h(z)|\bigr)\le
\frac{A_{\e}}{\e (r)} T\bigl((1+\varepsilon (r) )\cdot r\, ; f\bigr)
+B_{\e}\, , \quad r=|z|>0 \, ,
\end{equation}
holds. Here $A_{\e}$ and $B_{\e}$ are constants depending on the function
$\e (r)$ (but independent of $r$).
\end{Theorem}

In particular for  \, $1+\e (r) \equiv B>1$ in \eqref{Te} the last
Theorem \ref{Khab01} implies the Rubel--Taylor--Miles Theorem \ref{RTM},
and for \, $\e (r) =\epsilon \bigm/(1+r) $ in \eqref{Te} with a constant
$\epsilon >0$ we obtain the sharpening
of the one-dimensional variant of Skoda Theorem  \ref{Sko}
with $(1+ r)$ in place of $(1+r)^3$ in \eqref{Skoest}.
Besides the Theorem \ref{Khab01} can be  essentially stronger
than the Rubel--Taylor--Miles  Theorem or the one-dimensional Skoda Theorem
if  the Nevanlinna characteristic $T(r; f)$ has
rapid growth and the decreasing function $\e (r)$  is chosen in conformity
with the growth of $T(r; f)$.

For $n>1$  the first result on the representation problem
for a special growth functions $\lambda (r)=r^{\rho }$
(that solves the problem of L.~A.~Rubel \cite[Problem 8]{Problems} (1968))
is apparently due to W.~Stoll  \cite[Proposition 6.1]{Stoll} (1968).

Let $\lambda (\vec{r}\, )=\lambda (\Fl{1}{n}{r})$,
$\, \vec{r}\,  \in \R_+^n$, be a positive
continuous function which is nondecreasing in each variable,
and $T(\vec{r}\, , f)\leq \lambda (\vec{r}\, )$, $\, \vec{r}\, \in
\R_+^n$. B.~A.~Taylor  posed the problem\footnote{In the original,
``{\it When is $\Lambda$ the field of quotients of
${\Lambda}_{E}\,$\/}?''} \cite[p.~470]{Taylor} (1968):
\begin{description}
\item[{\bf (TP)}]
{\it when can a meromorphic function $f$ be represented as the quotient
$f=g/h$ of entire functions $g$ and $h$ such that
\begin{equation}
\label{estim2}
\log \bigl(|g(z)|+|h(z)|\bigr)\le A\lambda (B\vec{r}\, )+C \, ,
\; r_j=|z_j|\ge 1 \, , \; 1\leq j \leq n\, ,
\end{equation}
where $A$, $B$ and $C$ are constants ?\/}
\end{description}

He proved the following Theorem  \cite[Theorem]{Taylor} (1968):
\begin{Theorem}[B.~A.~Taylor]
\label{Tay}
If $\lambda (\vec{r}\, )$ is
slowly increasing in each variable, in the sense that
\begin{equation*}
\lambda (r_1, \ldots , 2r_j, \ldots , r_n)\le
A_j\lambda (\Fl{1}{n}{r})
\end{equation*}
for some constant $A_j>0$, $1\le j \le n$, and
$T(\vec{r}\, ; f)\le \const \cdot \lambda (\vec{r}\, )$
for all $r_j \ge 1$, $1\le j \le n$, then
there are entire functions $g$ and $h$ such that  $f=g/h$ and
\eqref{estim2} holds for all $r_j \ge 1$.
\end{Theorem}

Further progress for meromorphic functions $f$ in several variables and
for the characteristic $T(r; f)$, i.e. $k=1$,
under  special conditions on the behavior of $\lambda$,
 was made in \cite{KujalaB}, \cite[Propositions 7.3, 9.10]{Kujala} (1969-71):
\begin{Theorem}[R.~O.~Kujala]
\label{Kuj}
If there are constants $a_0$ $($resp., a vanishing function $a_0=a_0(r)$,
$r\to +\infty$$)$,
$a$, $b$ and $R$ in $\R_+$  and $p_0$ in $\N$ such that
\begin{equation*}
\int\limits_s^r \lambda(t)t^{-p-1}\, dt\le a_0\lambda (br)r^{-p}
+a\lambda (bs)s^{-p}
\end{equation*}
whenever $r\ge s >R$ and $p_0\le p$ in $\N$  $($in particular,
if $\lambda$ satisfies the slow growth condition \eqref{slow1}\/$)$, then
under the condition  \eqref{ballestm}
the meromorphic function  $f$ in $\bC^n$ can be represented as the quotient
$f=g/h$ of entire functions $g$ and $h$ such that
\eqref{balleste} $($recp., with a vanishing function $A=A(r)$, $r\to
+\infty$$)$ is valid.
\end{Theorem}

Except the Theorems \ref{LN}, \ref{Li}, \ref{Gol}, \ref{Khab01} and the
mentioned result for $\lambda (r)=r^{\rho}$ from W.~Stoll \cite{Stoll},
all these results were obtained by the Fourier series method
(see  \cite{Kondratyuk}, \cite{Kujala}, \cite{Miles},
\cite{Noverraz}, \cite{Taylor-Rubel}--\cite{Rubel}, \cite{Taylor}).
This method was used first by N.~I.~Akhiezer in \cite{Akhiezer} (1927)
for the new proof  the Lindel\"of Theorem and
later  by L.~A.~Rubel in \cite{Rubel1} (1961) (see in this
connection also \cite[P.~85--88]{GO}).

Using the $\bar \partial$-problem method, H.~Skoda
\cite{Skoda}--\cite{Skoda72} (1971-72) obtained  the following
results (see also \cite[Theorem 9.12]{Stoll85}):
\begin{Theorem}[H.~Skoda]
\label{Sko}
For every  constant $\epsilon >0$  entire functions $g$ and $h$
in representation $f=g/h\in \Mer$ can be chosen to
satisfy the estimates
\begin{equation}
\label{Skoest}
\log \bigl(|g(z)|+|h(z)|\bigr)
\le C(\epsilon ,s)(1+r)^{4n-1}T(r+\epsilon ; f)
\end{equation}
or
\begin{equation}
\label{Skoest1}
\log \bigl(|g(z)|+|h(z)|\bigr)\le C(\epsilon , s)
\bigl(\log (1+r^2)\bigr)^2 T\bigl((1+\epsilon )r; f\bigr)
\end{equation}
for all $r=|z|\geq s>0$,
where $C(\epsilon , s)$ is a constant depending on $\epsilon $ and
$s$.
\end{Theorem}

It should be noted  \cite[p.~295]{Stoll85} that the pair $(g, h)$
in the Theorem \ref{Sko} can have a common divisor, can depend on
$\varepsilon$ but can not depend on $s$, and can be different in
\eqref{Skoest1} from the pair chosen in \eqref{Skoest}. The case
\eqref{Skoest} is good for rapid growth, the case \eqref{Skoest1}
is good for slow growth.

An analysis  of these  results indicates a definite
rift between the cases of slow and rapid growth of the majorants
$\lambda$ or of the characteristics $T(\cdot ; f)$, both in methods
(the Fourier series method in the Theorems \ref{RTM}, \ref{Tay}--\ref{Kuj}
and the $\bar \partial$-problem method in the Theorem \ref{Sko})
and  the estimates and conditions on $\lambda$ or $T(\cdot
; f)$. Our balayage (or sweeping out) method (1990-93) enables
 us to get rid of this rift and yields results in a complete form.

To be more specific, first this balayage method was applied
for a new short nonconstructive proof of Rubel-Taylor-Miles Theorem
\ref{RTM}  in \cite[\S~3]{Khabib91} (1991).
The improvement of this method for $n=1$ made it possible
recently to prove the Theorem \ref{Khab01}.
In our articles \cite{Khabib92}, \cite{Khabib93-I}
we established  for functions of several variables
 the following results that are, in a certain sense, extreme  relative  to
the Nevanlinna characteristics $T(r; f)$, $r\geq 0$, and \eqref{Tcirc}.

\begin{Theorem}[B.~N.~Khabibullin]
\label{repres}
Let $f\in \Mer$, $n\ge 1$. Then
\begin{description}
\item[(T)]
{\rm {\cite[Theorem 1.3]{Khabib92}}, {\cite[Theorem 5]{Khabib93-I}}
(1992-93)}

For each constant $\e >0$  the function $f$ can be represented
as the quotient
\begin{equation} \label{gh}
f=\frac{g}{h}\, , \qquad g \text{ and } h \text{ are entire functions},
\end{equation}
such that
\begin{equation}
\label{estT}
\log \bigl(|g(z)|+|h(z)|\bigr)\le A_{\e}T\bigl((1+\e )r; f\bigr)+C_{\e}\, ,
\quad     r=|z|>0;
\end{equation}

\item[(Tc)]
{\rm {\cite[Theorem 4]{Khabib93-I}} (1993)}

Under the conditions  $T_f^c(z)\le \lambda (z)$,
$z\in \bC^n$, $f(0)=1$, where the function  $\lambda$  is circular
continuous positive increasing on all rays with origin at $0$
and satisfies the following two H\"ormander conditions\/
{\rm \cite{Hormander67}:}

\begin{description}
\item[(L)] \label{L}
$ \log \bigl(1+|z|\bigr)\leq O\bigl(\lambda (z)\bigr)\; \text{ as } \;
      |z|\to +\infty \,$,
and
\item[(H)] \label{Hor}  for any number
  $\epsilon >0$ there exist positive constans $c_1, c_2, c_3, c_4 $ such that
$$
|z-\zeta |\le \exp \bigl(-c_1\lambda (z)-c_2\bigr) \Longrightarrow \,
\lambda (\zeta)\le c_3\lambda\bigl((1+\epsilon )z\bigr)+c_4,
\; z\in \bC^n,
$$
or
\item[(\^H)]
for any $\epsilon  >0$ there exist positive constants $\sigma$, $c_1$, and
$c_2$ such that
\begin{equation}   \label{H2}
|z-\zeta |\le \sigma \, \Longrightarrow \,
\lambda (\zeta)\le c_1\lambda\bigl((1+\epsilon )z\bigr)+c_2,  \quad z\in \bC^n,
\end{equation}
\end{description}
for any $\e >0$  the function $f$ can be represented
as a quotient  \eqref{gh} such that
\begin{equation} \label{cTlambda}
\log \bigl(|g(z)|+|h(z)|\bigr)\le A_{\e} \lambda \bigl((1+\e ) z\bigr)
+ C_{\e}\; ,
\quad z\in {\bC}^n \; ,
\end{equation}
where under the condition \eqref{H2}
functions $g$ and $h$ can be chosen so that $g(0)=h(0)=1$.
\end{description}

Here $A_{\e}$ and $C_{\e}$ are constants depending on $\e$.
\end{Theorem}

Besides the fact that the extra conditions on the growth function
$\lambda$ was removed in  the Theorem \ref{repres} (Part {\bf (T)}),
 it also refines
the Theorems \ref{RTM}, \ref{Kuj} and \ref{Sko}. For example, the constant
$B$ in  \eqref{balleste} and \eqref{ballesteM}
is replaced by a constant $1+\e >1$
arbitrary close to $1$, and in \eqref{Skoest} and in \eqref{Skoest1}
it is possible to
remove the power factor before  $T(\cdot ; f)$
if in \eqref{Skoest} $r+\e$  is replaced by $(1+\e )r$, where $\e >0$.
Following Gol'dberg \cite{Gol'dberg}, we can show that we cannot set
$\e =0$ in general in \eqref{estT} and in \eqref{cTlambda}.
An analog of the Theorem \ref{repres} (Part {\bf (T)})
for $\delta$-subharmonic functions was established
by O.~V.~Veselovskaya \cite{Ves} (1984) by the Fourier series method.

If we have  some additional information on the
functions $g_f$ and $h_f$ in  the initial representation
\eqref{repr}, then these additional properties can be preserved
sometimes in the final representation \eqref{gh}
(see \cite[Theorems 1--3]{Khabib98} (1998)):

\begin{Theorem}[B.~N.~Khabibullin] If
it is known in addition to conditions of  Theorem \ref{repres} that for
$f\in \Mer$ with $f(0)=1$  there exists
the representation  \eqref{repr}, where $g_f$ and
$h_f$   are  \underline{bounded in some open set} $B\subset {\bC}^n$,
then in  {\rm {\bf (T)}}  $($recp., in {\rm {\bf (Tc)}}$)$
functions $g$ and $h$  can be chosen  so that
besides  \eqref{estT} $($recp., \eqref {cTlambda}$)$
for any $\gamma \geq 1$ the inequlity
$$
\log \bigl(|g(z)|+|h(z)|\bigr)\le C_{\gamma }\log \bigl(2+|z|\bigr)\; ,\quad
z\in B_{\gamma }\, ,
$$
holds, where
$$
B_{\gamma }=\{ z\in B : \; \dist (z,\partial B)\ge
\bigl(1+|z|\bigr)^{-\gamma}\}
\, ,
$$
$\dist (z,\partial B)$ is the distance from $z$ to the boundary
$\partial B$ of $B$,  and $C_{\gamma }$ is constant.

If $n=1$, then functions $g$ and $h$ in \eqref{gh} can be chosen
also so that $g$ and $h$ are bounded on the set
$\{z\in B : \dist (z,\partial B)\ge 1/{\gamma}\}$.
\end{Theorem}

For the general case of the Nevanlinna $\mathbf{\mbox{\L}}$-characteristic
we have at present  the following result \cite{Khabib00} (2000):
\begin{Theorem}[B.~N.~Khabibullin] \label{TL}
Let $f\in \Mer$ and $\bC^n$
is represented as \eqref{directsum}. Then,
with the preceding notations \eqref{uf}--\eqref{T},
for any constant $\e >0$,   there exists a representation
\eqref{gh} such that
\begin{equation} \label{TLest}
\log \bigl(|g(z)|+|h(z)|\bigr)
\le A_{\e} \text{\rm $T_{\mathbf{\mbox{\L}}}^+$}\bigl((1+\e )\vec r\, +
B_{\e}\cdot \vec 1 \, ; f\bigr)
+C_{\e}\log \bigl(2+| z|\bigr)
\end{equation}
for all $z\in \bC^n$, where  $A_{\e}, B_{\e}, C_{\e}$ are constants and
$\vec 1\, =(1, \ldots , 1) \in\R_+^k$, $\; T^+=\max \{ T, 0\}$.
\end{Theorem}

It is easy to see, that the Theorem \ref{TL} implies the Theorem
\ref{repres} (Part {\bf (T)}),
because always $\log \bigl(2+| z|\bigr)\leq O\bigl(T(r, f)\bigr)$,
$r=|z|$, $f\not\equiv 0, \infty$.
Also, if
\begin{equation} \label{log}
\log \bigl(2+|\vec r\,|\bigr)\leq O\bigl(\lambda (\vec{r}\,
)\bigr)\, , \quad  k=n, \quad |\vec r \, |=\sqrt{r_1^2+\cdots +r_k^2}\, ,
\end{equation}
then  the Theorem \ref{TL}  is the extension of Theorem \ref{Tay}.
It  solves the Taylor problem {\bf (TP)} with  minimal
estimate \eqref{log} from below for $\lambda (\vec{r}\, )$.

A general scheme of the solution of the representation
problem and some other problems (balayage method) was presented in
\cite{Khabib96}--\cite{Khabib97}, \cite{Khabib01I}.

\subsubsection*{Unsolved problems.}

\paragraph*{Problem 1.} Is it true,
that {\it the summand $C\log \bigl(2+| z|\bigr)$
in  \eqref{TLest} in the Theorem \ref{TL}
can be removed\/} ?

In any case it is true for $k=1$, and in the Theorem \ref{Tay} for $k=n$
if $\lambda (\vec r \, )$ is slowly increasing.

\paragraph*{Problem 2.}
{\it Extend the result of the Theorem \ref{Khab01}
to $\bC^n$, $n>1$, for classic Nevanlinna characteristic $T(r; f)$
and for Nevanlinna $\mbox{\L}$-characteristic.}

\paragraph*{Problem 3.}
Denote by   $\mathcal{PSH}_p$  the cone of all
{\it plurisubharmonic\/} function on $L_p$, $p=1, \ldots , k$,
where $L_p$ was defined in \eqref{directsum}.
Denote by  $\mathcal{M}_p^+$ the cone of all {\it positive Borel measures
$($Radon measures$)$ with compact support\/} in $L_p\, .$

We introduce in $\mathcal{M}_p^+$ the {\it partial order\/}
${\prec}$ by (see \cite[Ch.~XI]{Meyer})
$$
 (\; \nu \,{\prec} \,\mu\;  )\; \Longleftrightarrow \;
\Bigl(\;  \int v\, d\nu \leq \int v\, d\mu
\;  \text{ for all } \; v\in \mathcal{PSH}_p \; \Bigr)\, .
$$
If $\nu \prec \mu$, then the measure $\mu\in \mathcal{M}_p^+$,
is called a {\it balayage\/}
( also  {\it sweeping out\/} or {\it Jensen measure\/} \cite{Gamelin})
of the measure $\nu$ (with respect to the cone $\mathcal{PSH}_p\,$).

Further,  let
$\boldsymbol \mu =(\Fl{1}{k}{\mu})$
be a fixed simply ordered collection of measures, where  every measure
$\mu_p\in \mathcal{M}_p^+$ is  the balayage
of the Dirac measure $\delta_p \in \mathcal{M}_p^+$
at the point $0\in L_p$, i.e. $\int u \, d\delta_p =u(0)$.

Let $u_f$ be a plurisubharmonic function from \eqref{uf}.
A function
\begin{equation*}
T_{\mathbf{\mbox{\L}}, \boldsymbol \mu }(\vec r\, ; f)=
\int\limits_{S_1} \cdots \int\limits_{S_k} u_f(r_1\zeta_1, \ldots , r_k\zeta_k)
\, d\mu_1 ({\zeta}_1)\cdots d\mu_k ({\zeta}_k)\, ,
\quad  \vec r\, \in \R_+^k,
\end{equation*}
will be called the {\it Nevanlinna-Jensen\/}
$(\mathbf{\mbox{\L}}, \boldsymbol \mu )$-{\it characteristic\/} of  $f\in \Mer$.
{\it How can the Theorem \ref{TL} be extended  for the
Nevanlinna-Jensen
{\rm $(\mathbf{\mbox{\L}}, \boldsymbol \mu )$}-cha\-r\-a\-c\-t\-e\-r\-i\-s\-tic\/} ?

\section{The representation of meromorphic\\
functions with restrictions on the type\\
and on the circled indicator}

\setcounter{equation}{0}

In this section we consider the problem of the representation of a
meromorphic function $f$ in $\bC^n$ as a quotient \eqref{gh}
with  best possible estimates of the circled indicators and of the types
of the entire functions $g$ and $f$.

{\it The representation} (or {\it factorization}) {\it type
${\sigma}_n^*(f, \rho )$ of a meromorphic function $f$ on $\bC^n$ of  order
$\rho$} is defined to be the infimum of the numbers $\sigma$
for which there are entire functions $g$ and $h$ of
order $\rho$ with the type less
than $\sigma$  (see \cite{LG}) such that $f=g/h$.

Let
\begin{equation} \label{Prho}
P_1 (\rho )=\begin{cases}
\dfrac{\pi \rho}{\sin \pi \rho} \, ,
& \text{ if \; $0\leq \rho \leq 1/2$}\, , \\
\pi \rho \, , & \text{ if \; $\rho \, >\, 1/2$}\, .
\end{cases}
\end{equation}

For $n=1$ the first sharp result in this direction was established
apparently in our paper \cite[Theorem 2]{Khabib92t} (1992):

\begin{Theorem}[B.~N.~Khabibullin]  \label{Tn1}
Let $f$ be a meromorphic
function  in the complex plane $\bC$. If the type of the Nevanlinna
characteristic $T(r; f)$ is equal to $\sigma$ with order $\rho$, i.e.
\begin{equation} \label{Ttype}
\limsup_{r\to +\infty} r^{-\rho}T(r;f)=\sigma\, ,
\end{equation}
then the sharp estimate
$\;
\sigma \leq  {\sigma}_1^*(f, \rho )  \leq \sigma P_1 (\rho )
\;$
is valid.
\end{Theorem}

For $n>1$ we set
\begin{equation} \label{Pnrho}
P_n (\rho )=P_1 (\rho )\, \prod_{k=1}^{n-1} \Bigl(1+\frac{\rho }{2k}\Bigr)\, .
\end{equation}

For $n>1$ we have \cite[Theorem 3]{Khabib93-II} (1993):
\begin{Theorem}[B.~N.~Khabibullin] \label{Tn}
If the type of the Nevanlinna
characteristic $T(r; f)$ is equal to $\sigma$ with order $\rho$, i.e.
 \eqref{Ttype} holds,
then
$$
\sigma \leq  {\sigma}_n^*(f, \rho )  \leq
\sigma P_1 (\rho ) \, \prod_{k=1}^{n-1}b\Bigl(\frac{\rho}{2k}\Bigr)
\leq e^{n-1}P_n(\rho ) ,
$$
where
\begin{equation}   \label{Pb}
b(x)=\begin{cases}
ex \, ,                  & \text{ if \; $x\geq 1$}\, , \\
e^x \,  ,                & \text{ if \; $x< 1$}\, ,
\end{cases}
\end{equation}
and for $\rho \leq 1$ the sharp estimate

\begin{equation}   \label{Psh}
\sigma \leq  {\sigma}_n^*(f, \rho )  \leq \sigma P_n (\rho )
\end{equation}
is valid.
\end{Theorem}

Let $u$ be a plurisubharmonic function on $\bC^n$ of finite type with
order $\rho$. Define the function
$$
h_{c\rho}(z, u)=\limsup_{\xi \to \infty} |\xi |^{-\rho}u(\xi z),
\quad \xi \in \bC \, , \; z\in \bC^n\, .
$$
Its upper semicontinuous regularization
$$
h_{c\rho}^*(z, u)=\limsup_{\zeta \to z} h_{c\rho}(\zeta , u)
$$
is called \cite[Definition 1.29]{LG} the {\it circled indicator} of
growth of $u$ (with order $\rho$).
If we are dealing with circled indicator of growth  $h_{c\rho}^* (\zeta , u)$,
then we assume a priori that $u$ is of finite type with order $\rho$.
It is known that the circled indicator is complex-homogeneous of
degree $\rho$, i.e., $h_{c\rho}^*(\xi z, u)=|\xi |^{\rho}
h_{c\rho}^*(z, u)$, $\xi \in \bC$, is plurisubharmonic, and
$h_{c\rho}^*(z, u)\geq 0$ if $u\not\equiv -\infty$.
We remind that the circular Nevanlinna characteristic
$T_f^c(z)$ of $f\in \Mer$ in \eqref{Tcirc} is plurisubharmonic.
Therefore, circled
indicator of growth for $T_f^c$ is defined.
In terms of the circled indicator of growth $h_{c\rho}^* (\zeta , T_f^c)$
we have the following best possible result \cite[Theorem 5.1, Example
5.1]{Khabib93-II} (1993):

\begin{Theorem}[B.~N.~Khabibullin]  \label{Thcirc}
Suppose that $f\in \Mer$, $f(0)=1$,
and $k$ is a continuous circular function on $S^n \subset \bC^n$.
\begin{enumerate}
\item If
$\quad
h_{c\rho}^* (\zeta , T_f^c) <k(\zeta ), \quad \zeta \in S^n\, , \,
$
then there exist entire functions $g$ and $h$ such that $f=g/h$,
$g(0)=h(0)=1$, and
$$
\max\bigl\{ h_{c\rho}^* (\zeta , g), h_{c\rho}^* (\zeta , h)\bigr\}
<P_1(\rho )k(\zeta),
\quad \zeta \in S^n\, .
$$
The constant $P_1(\rho )$ cannot in general be diminished, not
even for any one of the functions $g$ and $h$.
\item If $f=g/h$ is a representation of $f$ as the quotient of entire
functions  $g$ and $h$, then
$$
\max\bigl\{ h_{c\rho}^* (\zeta , g), h_{c\rho}^* (\zeta , h)\bigr\}\geq
h_{c\rho}^* (\zeta , T_f^c)\, , \quad \zeta \in S^n\, ,
$$
and this estimate is sharp.
\end{enumerate}
\end{Theorem}

The main difficulty consists in  the obtaining  the upper estimates
in the Theorems    \ref{Tn1}--\ref{Thcirc}. The success here
was achieved with the use of the balayage method (see the preceding
section).

\subsubsection*{Unsolved problems.}
\paragraph*{Problem 1.}
{\it  Prove  the upper estimate in \eqref{Psh} for all $\rho$\/}.
See the commentary to the last Problem  in the next section.

\paragraph*{Problem 2.} {\it Extend the results of the Theorems
\ref{Tn}--\ref{Thcirc} to
the Nevanlinna {\rm $\mathbf{\mbox{\L}}$}-characteristic
and to the the Nevanlinna-Jensen {\rm $(\mathbf{\mbox{\L}},
\boldsymbol \mu )$}-characteristic\/}.

\section{Paley problem}

\setcounter{equation}{0}

The subject of this section is close to the themes of the previous
sections in both  estimates and methods of proofs.

For a   function $u$ in $\bC^n$ (recp., in $\R^m$)
with the range of values $[-\infty , +\infty]$
or $\R \cup \{ \infty\}$, we set
\begin{equation*}
M(r; u)=\max \bigl\{ u(z): |z|=r\} \, , \quad
u^+(z)=\max \bigl\{ u(z), 0 \bigr\}\, .
\end{equation*}

In 1932 R.~Paley supposed that for any entire function $f$ of order
$\rho \geq 0$ the inequality
\begin{equation} \label{Paleycon}
\liminf\limits_{r\to +\infty}
\frac{M\bigl(r; \log |f|\bigr)}{T(r; f)}
\leq P_1(\rho )
\end{equation}
holds (see the definition of  $P_1(\rho )$ in \eqref{Prho}).

In \eqref{Paleycon} the equality is achieved, in particular, for the
Mittag-Leffler's function (see, for example,  \cite[p.111]{GO}).
Earlier the inequality \eqref{Paleycon} was proved by G.~Valiron
\cite{Valiron} (1930) and by A.~Whalund \cite{Whalund} (1929)
for $\rho \leq 1/2$.

The Paley conjecture was conclusively proved by N.~V.~Govorov
\cite{Govorov} only in 1967. In 1968, V.~P.~Petrenko \cite{Petrenko} extended
this result to meromorphic functions of finite lower order.
Different proofs  was obtained later by I.~V.~Ostrovski\u{\i}
\cite{Ostrovsk} for meromorphic
functions and by M.~R.~Ess\'en \cite{Essen} for subharmonic functions
in the complex plane.
\begin{Theorem}[N.~V.~Govorov--V.~P.~Petrenko] \label{GP}
If $f$ is a meromorphic function in the complex plane $\bC$
of finite lower order $\lambda$, i.e.,
\begin{equation} \label{lowerorder}
\lambda =\liminf_{r\to +\infty}
\frac{\log T(r; f)}{\log |r|} \, < \, +\infty \, ,
\end{equation}
then the inequality
\begin{equation} \label{Paleylambda}
\liminf\limits_{r\to +\infty}
\frac{M\bigl(r; \log |f|\bigr)}{T(r; f)}
\leq P_1(\lambda )
\end{equation}
holds.
\end{Theorem}

The relative growth of $T(r; f)$ and $M\bigl(r; \log |f|\bigr)$
for entire and meromorphic functions of infinite order
in the complex plane has also been  considered by many authors.

For entire functions of infinite order in the complex plane
in \cite{Chuang}(1981), \cite{MSh}(1984) and finally in
\cite{DDL}(1990--1991) was proved the following result.

\begin{Theorem}\label{DDL}
{\bf (C.~T.~Chuang, I.~I.~Marchenko--A.~I.~Shcherba,\linebreak
C.~J.~Dai--D.~Drasin--B.~Q.~Li)}
Let  $f$ be an entire function of infinite order in the complex plane.
If $\psi (x)$ is  increasing and positive for $x\geq x_0>0$
and if
\begin{equation}\label{psi}
\int_{x_0}^{\infty}\frac{dx}{\psi (x)}<\infty \, ,
\end{equation}
then
\begin{equation}\label{Tpsi}
\liminf_{r\to +\infty}\frac{M\bigl(r; \log |f|\bigr)}
{T(r; f) \cdot \psi\bigl( \log T(r; f)\bigr)}=0
\end{equation}
and even
\begin{equation*}\label{oTpsi}
M\bigl(r; \log |f|\bigr)=
o\Bigr( T(r; f) \cdot \psi\bigl( \log T(r; f)\bigr)\Bigr)
\end{equation*}
as $r\to +\infty$  through a set of logarithmic density one.
\end{Theorem}
In \cite{MSh} and \cite{DDL} it is also shown
that the results are best possible in some sense.

For meromorphic  functions
in \cite{DDL} (1990--1991) the following was obtained.
\begin{Theorem}[C.~J.~Dai-D.~Drasin-B.~Q.~Li]\label{DDLm}
Let $f$ be a meromorphic fu\-n\-c\-tion
of infinite order in the complex plane.
If $\psi (x)$ is the same as in the Theorem \ref{DDL}
then
\begin{equation*}\label{oTmpsi}
M\bigl(r; \log |f|\bigr)=
o\Bigr( T(r; f) \cdot \psi\bigl( \log T(r; f)\bigr)
\cdot \log \psi\bigl( \log T(r; f)\bigr)\Bigr)
\end{equation*}
as $r\to +\infty$  through a set of logarithmic density one.
\end{Theorem}

A different approach has been taken in \cite{Bergweiler}, \cite{BB}
(1990--1994) where the characteristic $M\bigl(r; \log |f|\bigr)$  has been
compared  with the derivative of $T(r; f)$.

\begin{Theorem}[W.~Bergweiler-H.~Block]\label{BBm}
Let $f$ be a meromorphic fun\-c\-tion of infinite order in the complex
plane. Then
\begin{equation*}\label{oTgam}
\liminf_{r\to +\infty}
\frac{M\bigl(r; \log |f|\bigr)}{r T_{-}'(r; f)}\leq \pi \, ,
\end{equation*}
where $T_{-}'(r; \cdot )$ is the left-side derivative of $T$ of $r$.

Let $\psi (x)$ be positive and continuously
diffirentiable for $x\geq x_0 >0$ such that $\psi (x)/x$
is non-decreasing, $\psi (x)\leq \sqrt{\psi (x)}$, and
\eqref{psi} is satisfied. Then \eqref{Tpsi} holds.
\end{Theorem}

The articles of I.~I.~Marchenko \cite{March1}-\cite{March2} contains
much  information in particular on the growth of entire and
meromorphic functions of infinite order.

We don't know any results for entire and meromorphic functions
of infinite order in $\bC^n$,
$n>1$, wich are analogs of the Theorems \ref{DDL}--\ref{BBm}.

Let $u$ be a subharmonic function in $\R^m$, $m\geq 2$,
(resp., in $\bC^n$)
and let $|\mathcal S^{m-1}|$ be the area of the unit sphere
$\mathcal S^{m-1}$ in $\R^m$, $d{s}$
is the area element on the unit sphere $\mathcal S^{m-1}$.
\begin{equation*}
m_q (r; u)=\biggl(\;  \frac1{|\mathcal S^{m-1}|}\,
\int\limits_{|\mathcal S^{m-1}|} \bigl|u(rx)\bigr|^q \,
d{s}(x) \,\biggr)^{1/q}\, ,
\quad r>0, \; 1\leq q <+\infty.
\end{equation*}
Then Nevanlinna characteristic $T(r;u)$ and $M(r; u)$ are respectively
$$
 T(r; u)=m_1 (r; u^+), \quad   M(r; u^+)=m_{\infty} (r; u^+)=
\lim_{q\to +\infty}m_q (r; u^+) .
$$
In particular if $f$ is a entire function in $\bC^n$ then
$T(r; f)=T\bigl(r; \log |f|\bigr)$, $m=2n$.

An analogue of \eqref{Paleylambda} for subharmonic functions
of finite lower order $\lambda$ in $\R^m$, $m\geq 3$, was
obtained by B.~Dahlberg \cite[Theorem 1.2]{Dahlberg}.

To be more specific suppose $\lambda \in (0, +\infty )$ is given.
The Gegenbauer functions $C_{\lambda}^{\gamma}$ are given as the solutions
of the differential equation
$$
(1-x^2)\, \frac{d^2u}{dx^2}-(2\gamma +1)x\, \frac{du}{dx}+\lambda
(\lambda+2\gamma )\, u=0, \quad -1<x<1\, ,
$$
with the normalization
$$
\lim_{x\to 1-0}C_{\lambda}^{\gamma} (x)=C_{\lambda}^{\gamma} (1)=
\frac{\Gamma (\lambda+2\gamma )}{\Gamma (2\gamma )\Gamma (\lambda +1)}.
$$
Put
$\;
a_{\lambda}=\sup \bigl\{ t: C_{\lambda}^{\frac{m-2}{2}} (t)=0\bigr\}
\; $
and define the function $u_{\lambda}$ in $\R^m$, $m\geq 3$, by
\begin{equation*}
u_{\lambda}(x)=\begin{cases}
0 \quad &\text{if $x_1\leq a_{\lambda}r$} \\
r^{\lambda}C_{\lambda}^{\frac{m-2}{2}} (x_1/r)
&\text{if $x_1> a_{\lambda}r$}\, ,
\end{cases}
\end{equation*}
where $x=(x_1, \ldots , x_m)$ and $r=|x|$.

The function $u_{\lambda}$ is subharmonic in $\R^m$ and the lower
order of $u_{\lambda}$ is $\lambda$. Set
\begin{equation*}
c(\lambda , m)=
\liminf_{r\to +\infty} \,
\frac{M(r; u_{\lambda})}{T(r; u_{\lambda })}\, .
\end{equation*}
\begin{Theorem}[B.~Dahlberg] \label{Dahl}
Let $u$ be a subharmonic function in
$\R^m$, $m\geq 3$, of finite lower order $\lambda >0$.
Then we have that
\begin{equation*}
\liminf_{r\to +\infty} \,
\frac{M(r; u )}{T(r; u)}\, \leq \, c(\lambda ,
m)\, ,
\end{equation*}
and this estimate is  best possible.
\end{Theorem}

The generalization of the Theorems \ref{GP}, \ref{Dahl}
is also known.

\begin{Theorem}[M.~L.~Sodin {\rm \cite{Sodin} (1983)}] \label{ThSodin}
Let $u$ be a subharmonic function
of lower order $\lambda$ in the complex plane. Then
\begin{equation} \label{mqSodin}
\liminf_{r\to +\infty}\frac{m_q (r; u^+)}{T(r; u)}\leq
m_q(S_{\lambda}), \quad 1 < q\leq +\infty,
\end{equation}
where $m_q(S_{\lambda})$ is the Lebesque mean of order $q$,
$1<q<+\infty$, of the function
\begin{equation*}
S_{\lambda}(\varphi )=\begin{cases} \pi \lambda \cos \lambda \varphi
\quad &\text{if $\quad |\varphi|\leq  \dfrac{\pi}{2\lambda}$},\\
0 \quad &\text{if $\quad \dfrac{\pi}{2\lambda}<|\varphi|\leq \pi $},
\end{cases}
\end{equation*}
for $\lambda \geq 1/2$ and
$$
S_{\lambda}(\varphi )= \frac{\pi \lambda \cos \lambda \varphi}{\sin
\pi \lambda}\, , \qquad |\varphi |\leq \pi,
$$
for $\lambda <1/2$, $ \; m_{\infty}(S_{\lambda})=
\max \bigl\{ S_{\lambda}(\varphi ): |\varphi |\leq \pi \bigr\}$.
The inequality \eqref{mqSodin} is sharp.
\end{Theorem}

This result was extended to $m\geq 3$ in \cite[Theorem 1]{KonTarVas}
(1995).

Set
\begin{equation*}
Q_{\lambda}^{\frac{m-2}{2}}(\theta )=\begin{cases}
A(\lambda, m)\, C_{\lambda}^{\frac{m-2}{2}}(\cos \theta )
\quad &\text{if $\quad 0\leq \theta \leq  \alpha_{\lambda}$},\\
0 \quad &\text{if $\quad \alpha_{\lambda}<\theta \leq \pi $},
\end{cases}
\end{equation*}
where
\begin{align*}
\alpha_{\lambda}&=\min \bigl\{\theta \in (0, \pi ):
C_{\lambda}^{\frac{m-2}{2}}(\cos \theta )=0 \bigr\}, \\
A(\lambda, m)&=\biggl( \; \frac{|\mathcal S^{m-2}|}{|\mathcal S^{m-1}|}
\, \int\limits_0^{\alpha_{\lambda}} \,
C_{\lambda}^{\frac{m-2}{2}}(\cos \theta ) {\sin}^{m-2}\theta \, d\theta
\biggr)^{-1}\, .
\end{align*}

\begin{Theorem}[A.~A.~Kondratyuk--S.~I.~Tarasyuk--Ya.~V.~Vasyl'kiv]
Let $u$ be a subharmonic function
of lower order $\lambda >0$ in $\R^m$, $m\geq 3$. Then
for every $q$, $1<q\leq +\infty$,  the inequality
\begin{equation} \label{mqKTV}
\liminf_{r\to +\infty}\frac{m_q (r; u^+)}{T(r; u)}\leq M_q({\lambda}, m)
\end{equation}
is true, where $M_q({\lambda}, m)=m_q
\bigl(Q_{\lambda}^{\frac{m-2}{2}}\, \bigr)$.
There exists a subharmonic function $u$ of lower order $\lambda >0$ in $\R^m$
for which in \eqref{mqKTV}  the equality is achieved.
\end{Theorem}

Besides, from W.~Hayman result \cite{Hayman} it follows that for
subharmonic functions of order $\lambda =0$ the inequlity
$$
\liminf_{r\to +\infty}\frac{m_q (r; u^+)}{T(r; u)}\leq 1, \quad
1<q\leq +\infty \, ,
$$
holds.

In \cite{Khabib95} (1995) we extended \eqref{Paleylambda} to
meromorphic functions $f$ in $\bC^n$, $n>1$. Instead of $M(r; \log |f|)$
the corresponding formulation involves in this case the family of
characteristics $M(r; \log |f_{\zeta}|)$, $\zeta \in S^n$, of the
slice functions $f_{\zeta}(w)=f(w\zeta )$, $w\in \bC$. This is
understanble, for in the case of a meromorphic function the quantity
$M(r; \log |f|)$ can be identically equal to $\infty$ starting from
some value of $r$, whereas the Nevanlinna characteristics of $f$ can
have very slow growth (one example is the function $f(z_1,
z_2)=(1+z_1)/(1+z_2)$ in $\bC^2$).
\begin{Theorem}[B.~N.~Khabibullin] \label{Khabmer}
Let $f$ be a meromorphic function in
$\bC^n$ of finite lower order $\lambda$, i.e., \eqref{lowerorder} holds.
Then
\begin{equation} \label{Paleylan}
\liminf\limits_{r\to +\infty}
\frac{M\bigl(r; \log |f_{\zeta}|\bigr)}{T(r; f)}
\leq P_n(\lambda )\, , \quad \zeta \in S^n,
\end{equation}
where $P_n(\lambda )$ was defined in \eqref{Pnrho},
and this estimate is  best possible.
\end{Theorem}

The estimate $\eqref{Paleylan}$ was new also for entire functions and
was not implied in general by the above result of
B.~Dahlberg for subharmonic functions in $\R^m$, $m\geq 3$.

On the other hand, one can  pose the problem of finding a complete
and best possible analogue of \eqref{Paleylambda} for plurisubharmonic
functions and entire functions of several variables. Such result was
obtained in our paper \cite[Theorem 1]{Khabib99} (1999):

\begin{Theorem}[B.~N.~Khabibullin] \label{Khabpl}
Let $u$ be a plurisubharmonic function in
$\bC^n$ of finite lower order $\lambda$.
Then
\begin{equation} \label{Paleypsh}
\liminf\limits_{r\to +\infty}\frac{M(r; u)}{T(r; u)}\leq P_n(\lambda )
\end{equation}
for $\lambda \leq 1$  and this estimate is  best possible. For $\lambda >1$,
\begin{equation}  \label{Paleymore1}
\liminf\limits_{r\to +\infty}\frac{M(r; u)}{T(r; u)}\leq
P_1 (\rho ) \, \prod_{k=1}^{n-1}b\Bigl(\frac{\rho}{2k}\Bigr)
\leq e^{n-1}P_n(\rho ) ,
\end{equation}
where $b(x)$ was defined in \eqref{Pb}.
\end{Theorem}

More can be said about the definitive character of \eqref{Paleypsh}.
For each $\rho \geq 0$ there exist in $\bC^n$ entire functions of order
$\rho$ and normal type such that
$$
\lim_{r\to +\infty} \frac{M\bigl( r; \log |f|\bigr)}{T\bigl( r; \log
|f|\bigr)} = P_n(\rho ) \, .
$$

Comparing with the Dahlberg's Theorem  \ref{Dahl}, we can
draw the following  conclusions:
\begin{enumerate}[(a)]
\item  for $\lambda =0$ or $\lambda =1$ the result of the Theorem
\ref{Khabpl} is a consequence of Dahlberg's Theorem;
\item for $0<\lambda <1$ and $n=2$ Theorem \ref{Khabpl} does not follow
from  Dahlberg's Theorem;
\item for $\lambda \geq \pi e-1$ and $n=2$
even the estimate \eqref{Paleymore1} in Theorem \ref{Khabpl} does not
follow from Dahlberg's Theorem;
\item as $\lambda \to +\infty$, for $n=2$ the estimate \eqref{Paleymore1}
is better by order $O(1/\lambda )$ than the one following from
Dahlberg's Theorem.
\end{enumerate}

\subsubsection*{Unsolved problems.}

\paragraph*{Problem 1.}
{\it Extend the results of the Theorems \ref{DDL}--\ref{BBm}
for entire and meromorphic functions of several variables of infinite order.}

\paragraph*{Problem 2.}
{\it Extend the results of the Theorems
\ref{Khabmer}--\ref{Khabpl} to
the Nevanlinna {\rm $\mathbf{\mbox{\L}}$}-characteristic
and to the the Nevanlinna-Jensen {\rm $(\mathbf{\mbox{\L}},
\boldsymbol \mu )$}-characteristic\/}.

\paragraph*{Problem 3.}
{\it Extend the results of the  Theorem \ref{Khabpl} for
the $m_q(r; \cdot)$ instead of $M(r; \cdot )$.}

\paragraph*{Problem 4.}
{\it  Prove  the upper estimate in \eqref{Paleypsh} for all $\lambda$\/}.

\subparagraph*{Commentary.} In order to  prove the upper estimates
in \eqref{Psh} for all $\rho$ and in \eqref{Paleypsh}
for all $\lambda$, it is sufficient  to confirm the following:
\begin{Hypothesis}
Let $S$ be a  nonnegative increasing function on $\R_+$,
$S(0)=0$, and the function $S(t)$ is convex with respect to
$\log t$, i.e., $S(e^x )$ is convex on $[-\infty , +\infty )$.
Further, let $\lambda \geq 1/2$,  $n\in \N$, $n\geq 2$.
If
\begin{equation} \label{Hyppremise}
\int\limits_0^1 S(tx)(1-x^2)^{n-2}x\, dx \leq t^{\lambda}\, ,
\quad 0\leq t < +\infty \, ,
\end{equation}
then
\begin{equation} \label{Hypconcl}
\int\limits_0^{+\infty}
S(t)\frac{t^{2\lambda-1}}{(1+t^{2\lambda})^2}\, dt\leq \frac{\pi
(n-1)}{2\lambda}\prod_{k=1}^{n-1} \Bigl(1+\frac{\lambda}{2k}\Bigr)\, .
\end{equation}
\end{Hypothesis}

This Hypotesis is true if $\lambda \leq 1$. When
$$
S(t)=2(n-1)\prod_{k=1}^{n-1} \Bigl(1+\frac{\lambda}{2k}\Bigr)
\, t^{\lambda}\, , \quad \lambda \geq \frac12 \, ,
$$
we have the equalities in  \eqref{Hyppremise}  and in \eqref{Hypconcl}.

It is not known even, whether the formulated Hypothesis
is true for $n=2$ and $\lambda >1$ ?


\bigskip\par\noindent{\slshape
Department of Math., Bashkir State University,
Frunze St., 32, Ufa, Bashkortostan, 450074, Russia}

\bigskip\par\noindent{\slshape
Institute  of Math. of Ufa Scientific Center, Chernishevski\u{\i} St., 112,
Ufa, Bashkortostan, 450025, Russia}

\bigskip\par\noindent {\slshape E-mail: algeom$\, @\,$bsu.bashedu.ru
\hskip0.5cm and(or) \hskip0.5cm khabib-bulat$\, @\,$mail.ru}

\end{document}